\documentclass[12pt]{article}

\usepackage{graphicx,tikz,color,url}
\usepackage{amsmath,amssymb,latexsym}
\usepackage{MnSymbol}

\usepackage{tkz-graph}
\usepackage{tkz-berge}

\newtheorem{theorem}{Theorem}[section]

\newtheorem{proposition}[theorem]{Proposition}
\newtheorem{observation}[theorem]{Observation}

\newtheorem{corollary}[theorem]{Corollary}

\newcommand{\proof}{\noindent{\bf Proof.\ }}
\newcommand{\qed}{\hfill $\square$ \medskip}
\newcommand{\cp}{\,\square\,}

\newcommand{\diss}{{\rm diss}}
\newcommand{\diam}{{\rm diam}}
\newcommand{\gp}{{\rm gp}}
\newcommand{\sr}{{}_{\rm SR}}

\begin{document}

\title{General $d$-position sets}

\author{
Sandi Klav\v zar $^{a,b,c}$
\and
Douglas F. Rall $^{d}$
\and
Ismael G. Yero $^e$
}

\date{}

\maketitle

\begin{center}
$^a$ Faculty of Mathematics and Physics, University of Ljubljana, Slovenia \\
\medskip

$^b$ Faculty of Natural Sciences and Mathematics, University of Maribor, Slovenia \\
\medskip

$^{c}$ Institute of Mathematics, Physics and Mechanics, Ljubljana, Slovenia \\
{\tt sandi.klavzar@fmf.uni-lj.si}
\medskip

$^d$ Department of Mathematics, Furman University, Greenville, SC, USA\\
{\tt doug.rall@furman.edu}
\medskip

$^{e}$ Departamento de Matem\'aticas, Universidad de C\'adiz, Algeciras, Spain \\
{\tt ismael.gonzalez@uca.es}
\end{center}

\begin{abstract}
The general $d$-position number ${\rm gp}_d(G)$ of a graph $G$ is the cardinality of a largest set $S$ for which no three distinct vertices from $S$ lie on a common geodesic of length at most $d$. This new graph parameter generalizes the well studied general position number. We first give some results concerning the monotonic behavior of ${\rm gp}_d(G)$ with respect to the suitable values of $d$. We show that the decision problem concerning finding ${\rm gp}_d(G)$ is NP-complete for any value of $d$. The value of ${\rm gp}_d(G)$ when $G$ is a path or a cycle is computed and a structural characterization of general $d$-position sets is shown. Moreover, we present some relationships with other topics including strong resolving graphs and dissociation sets. We finish our exposition by proving that ${\rm gp}_d(G)$ is infinite whenever $G$ is an infinite graph and $d$ is a finite integer.
\end{abstract}


\noindent
{\bf AMS Subj.\ Class.\ (2010)}: 05C12, 05C63, 05C69

\section{Introduction}
\label{sec:intro}

A {\em general position set} of a graph $G$ is a set of vertices $S\subseteq V(G)$ such that no three vertices from $S$ lie on a common shortest path of $G$. The order of  a largest general position set, shortly called a {\em gp-set}, is the {\em general position number} $\gp(G)$ of $G$ (also written gp-number). This concept was recently and independently introduced in~\cite{manuel-2018a, ullas-2016}. We should mention though that  the same concept was studied on hypercubes already in 1995 by K\"orner~\cite{korner-1995}. Following~\cite{manuel-2018a} and its notation and terminology, the concept received a lot of attention, see the series of papers~\cite{ghorbani-2019, klavzar-2019+, klavzar-2020+, klavzar-2019, manuel-2018b, neethu-2020+, patkos-2019+}. In particular, in~\cite{neethu-2020+} the general position problem was studied on complementary prisms. In order to characterize an extremal case for the general position number of these graphs, the concept of general $3$-position was introduced as an essential ingredient of the characterization. In this paper we extend this idea as follows.

Let $d\in {\mathbb N}$ and let $G$ be a (connected) graph. Then $S\subseteq V(G)$ is a {\em general $d$-position set} if the following holds:
\begin{equation}
\label{eq:d-general-position}
\{u,v,w\}\in \binom{S}{3}, v\in I_G(u,w) \Rightarrow d_G(u,w) > d\,,
\end{equation}
where $d_G(u,w)$ denotes the shortest-path distance in $G$ between $u$ and $w$, and $I_G(u,w) = \{x\in V(G):\ d_G(u,w) = d_G(u,x) + d_G(x,w)\}$ is the {\em interval} between $u$ and $w$. In words, $S$ is a general $d$-position set if no three different vertices from $S$ lie on a common geodesic of length at most $d$. We will say that vertices $u,v,w$ that fulfill condition~\eqref{eq:d-general-position} lie in {\em general $d$-position}. The cardinality of a largest general $d$-position set in a graph $G$ is the {\em general $d$-position number} of $G$ and is denoted by $\gp_d(G)$.

We proceed as follows. In the rest of this section we recall needed definitions and state some basic facts and results on the general $d$-position number. Then, in Section~\ref{sec:inequality-chain}, we demonstrate that in the inequality chain $\gp_{\diam(G)}(G) \le \gp_{\diam(G)-1}(G) \le \cdots \le \gp_{2}(G)$ all kinds of equality and strict inequality cases are possible. Using one of the corresponding constructions we also prove that the problem of determining the $\gp_d$ number is NP-complete. In Section~\ref{sec:paths-and-cycles} we determine the $\gp_d$ number of paths and cycles and give a general upper bound on the $\gp_d$ number in term of the diameter of a given graph. In the subsequent section we prove a structural characterization of general $d$-position sets. In Section~\ref{sec:connections} we report on the connections between  general $d$-position sets and two well-established concepts, the dissociation number and strong resolving graphs. In the concluding section we consider the $\gp_d$ number of infinite graphs and pose several open questions.

\subsection{Preliminaries}

For a positive integer $k$ we will use the notation $[k] = \{1,\ldots, k\}$. The \emph{clique number} and the \emph{independence number} of $G$ are denoted by $\omega(G)$ and $\alpha(G)$. If $S\subseteq V(G)$, then the subgraph of $G$ induced by $S$ is denoted by $\langle S\rangle$ and $\binom{S}{k}$ denotes the set of all subsets of $S$ having cardinality $k$. A subgraph $H$ of a graph $G$ is {\em isometric} if $d_H(u,v) = d_G(u,v)$ holds for all $u,v\in V(H)$. If $H_1$ and $H_2$ are subgraphs of $G$, then the distance $d_G(H_1,H_2)$ between $H_1$ and $H_2$ is defined as $\min\{d_G(h_1, h_2):\ h_1\in V(H_1), h_2\in V(H_2)\}$. In particular, if $H_1$ is the one vertex graph with $u$ being its unique vertex, then we will write $d_G(u,H_2)$ for $d_G(H_1,H_2)$. We say that the subgraphs $H_1$ and $H_2$ are {\em parallel}, denoted by $H_1\parallel H_2$, if for every pair of vertices $h_1\in V(H_1)$ and $h_2\in V(H_2)$ we have $d_G(h_1, h_2) = d_G(H_1, H_2)$. If $H_1$ and $H_2$ are not parallel, we will write $H_1\nparallel H_2$. The \emph{open neighborhood} and the \emph{closed neighborhood} of a vertex $v$ of $G$ will be denoted by $N_G(x)$ and $N_G[x]$, respectively. Vertices $x$ and $y$ of $G$ are {\em true twins} if $N_G[u] = N_G[v]$. We may omit the subscript $G$ in the above definitions if the graph $G$ is clear from the context.

Clearly, if $d = 1$, then every subset of vertices of $G$ is a general $1$-position set, and if $d\ge \diam(G)$, then $S$ is a general $d$-position set if and only if $S$ is a general position set. Moreover, note that
\begin{equation}
\label{eq:chain}
\gp_{\diam(G)}(G) \le \gp_{\diam(G)-1}(G) \le \cdots \le \gp_{2}(G)\,.
\end{equation}

We conclude the preliminaries with the following useful property.

\begin{proposition} \label{prop:isometric}
Let $G$ be a graph and let $2\le d\le \diam(G)-1$ be a positive integer. If $H_1,\dots,H_r$  are isometric subgraphs of $G$ such that $d_G(H_i,H_j) \ge d$
for $i\ne j$, then $\gp_d(G) \ge \sum_{i=1}^{r}\gp_d(H_i)$.
\end{proposition}

\proof
For each $i \in [r]$, let $S_i$ be a general $d$-position set of $H_i$ such that $|S_i|=\gp_d(H_i)$.  We claim that $S=\bigcup_{i=1}^r S_i$ is a general $d$-position set of $G$. Suppose $\{x,y,z\} \in \binom {S}{3}$ such that $y \in I_G(x,z)$ and $d_G(x,z)\le d$. That is, there exists a shortest $xz$-path of length at most $d$ in $G$ that contains $y$. Since $d_G(u,v)\ge d$ for any two vertices $u\in V(H_i)$ and $v\in V(H_j)$ with $i\ne j$, there exists $k\in [r]$ such that $\{x,y,z\}\subseteq V(H_k)$. Now, since $H_k$ is an isometric subgraph of $G$, it follows that  $d_{H_k}(x,y)=d_G(x,y)$, $d_{H_k}(y,z)=d_G(y,z)$ and $d_{H_k}(x,z)=d_G(x,z)$.  This implies that there is a $xz$-geodesic in $H_k$ that contains $y$.  Hence, $y \in I_{H_k}(x,z)$, and since $S_k$ is a general $d$-position set of $H_k$, we infer that
$d_G(x,z)=d_{H_k}(x,z) > d$, which is a contradiction.  Therefore, $S$ is a general $d$-position set of $G$, and it follows that $\gp_d(G) \ge \sum_{i=1}^{r}\gp_d(H_i)$.
\qed

\section{On the inequality chain~\eqref{eq:chain} and computational complexity}
\label{sec:inequality-chain}

In this section we investigate the inequality chain~\eqref{eq:chain} by constructing different classes of graphs which demonstrate that all kinds of equality and strict inequality cases can happen. We conclude the section by  applying one of these constructions to prove that the {\sc General $d$-Position Problem} is NP-complete.

\medskip\noindent
{\bf Equality in~\eqref{eq:chain} simultaneously}. \\
For $n \ge 2$, let $S$ be a star with center $x$ and leaves $u_1,\ldots,u_n,v_1,\ldots,v_n$.  Construct a graph $G_n$ of order $2n+3$ by taking the disjoint union of  $S$ and an independent set of vertices  $\{u,v\}$  together with the set of edges $\{uu_i,vv_i:\, i \in [n]\}$.  The diameter of $G_n$ is $4$, and we have $\gp_4(G_n)=\gp_3(G_n)=\gp_2(G_n)=2n$.

\medskip\noindent
{\bf Equality in~\eqref{eq:chain} simultaneously again}. \\
Let $T_r$, $r\ge 2$, be the tree obtained from $P_{r+1}$ by attaching two leaves to each of its internal vertices. Then we claim that
$$\gp_{r}(T_r) = \gp_{r-1}(T_r) = \cdots = \gp_{2}(T_r)\,.$$
Indeed, first note that $\diam(T_r) = r$. Since the gp-number of a tree is the number of its leaves (cf.~\cite[Corollary 3.7]{manuel-2018a}), we have $\gp_r(T_r) = 2r$. Let next $S$ be a general 2-position set. If $u$ is a vertex of $T_r$ adjacent to exactly two leaves, say $v$ and $w$, then $|S\cap \{u,v,w\}| \le 2$. Moreover, if $u$ is a vertex of $T_r$ adjacent to exactly three leaves, say $v$, $w$, and $z$, then $|S\cap \{u,v,w,z\}| \le 3$. It follows that $\gp_2(T_r) \le 2(r-3) + 2\cdot 3 = 2r$. In conclusion, $2r = \gp_{r}(T_r) \le \gp_{r-1}(T_r) \le \cdots \le \gp_{2}(T_r) \le 2r$, hence  equality holds throughout.

\medskip\noindent
{\bf Strict inequality in~\eqref{eq:chain} in exactly one case}. \\
Let $k, \ell \ge 4$ and let $G_{k,\ell}$ be a graph defined as follows. Its vertex set is
$$V(G_{k,\ell}) =  \bigcup_{j=1}^k \{u_j, w_j, x_{j,1},\ldots,  x_{j,\ell}\} \cup \{x_{k+1,1} \}\,.$$
For $j\in [k]$, each of the vertices $u_j$ and $w_j$ is adjacent to $x_{j,1},\ldots,  x_{j,\ell}$ and to $x_{j+1,1}$. There are no other edges in  $G_{k,\ell}$. Note that $|V(G_{k,\ell})| = k(\ell + 2) + 1$ and that $\diam(G_{k,\ell}) = 2k$.

It is straightforward to see that the set $X = \bigcup_{j=1}^k \{x_{j,1},\ldots,  x_{j,\ell}\} \cup \{x_{k+1,1} \}$ is a largest independent set of $G_{k,\ell}$. Moreover, $X$ is also a largest general $2$-position set and a largest general $3$-position set. Furthermore, it is not difficult to infer that the set $X\setminus \{x_{2,1},\ldots,  x_{k,1}\}$ is a largest general $d$-position set for each $d\in \{4, \ldots, 2k\}$. In conclusion,
$$\gp_{2k}(G_{k,\ell}) = \gp_{2k-1}(G_{k,\ell}) = \cdots = \gp_{4}(G_{k,\ell}) < \gp_{3}(G_{k,\ell}) = \gp_{2}(G_{k,\ell}) = \alpha(G)\,.$$

\begin{figure}[h]
\begin{tikzpicture}[scale=.5, transform shape]
\GraphInit[vstyle=Hasse]
  \SetUpEdge[color=black,lw=0.5pt]
  \begin{scope}[rotate=281.5]
    \grCycle[prefix=a,RA=3]{16}
    \EdgeInGraphMod{a}{16}{1}
    \EdgeInGraphMod{a}{16}{2}
    \EdgeInGraphMod{a}{16}{3}
    \EdgeInGraphMod{a}{16}{4}
    \EdgeInGraphMod{a}{16}{5}
    \EdgeInGraphMod{a}{16}{6}
    \EdgeInGraphMod{a}{16}{7}
    \EdgeInGraphMod{a}{16}{8}
  \end{scope}
  \node [scale=1.8] at (-2.5,-2.8) {$A$};

\draw[dotted, thick](3.2,-2.7)--(6.5,-2.7);
\draw[dotted, thick](3.8,0)--(6,0);
\draw[dotted, thick](3.2,-2.7)--(6,0);
\draw[dotted, thick](3.8,0)--(6.5,-2.7);
\draw[dotted, thick](3.2,2.7)--(6.5,2.7);
\draw[dotted, thick](3.8,0)--(6.5,2.7);
\draw[dotted, thick](3.2,2.7)--(6,0);
\draw[dotted, thick](3.2,-2.7)--(6.5,2.7);
\draw[dotted, thick](3.2,2.7)--(6.5,-2.7);
\end{tikzpicture}
\begin{tikzpicture}[scale=.5, transform shape]
\GraphInit[vstyle=Hasse]
  \SetUpEdge[color=black,lw=0.5pt]
  \begin{scope}[rotate=281.5]
    \grCycle[prefix=a,RA=3]{16}
    \EdgeInGraphMod{a}{16}{1}
    \EdgeInGraphMod{a}{16}{2}
    \EdgeInGraphMod{a}{16}{3}
    \EdgeInGraphMod{a}{16}{4}
    \EdgeInGraphMod{a}{16}{5}
    \EdgeInGraphMod{a}{16}{6}
    \EdgeInGraphMod{a}{16}{7}
    \EdgeInGraphMod{a}{16}{8}
  \end{scope}
  \node [scale=1.8] at (-2.5,-2.8) {$B$};
  \begin{scope}[rotate=0]
    \grEmptyCycle[prefix=b,RA=7]{1}
  \end{scope}
  \begin{scope}[rotate=0]
    \grEmptyCycle[prefix=c,RA=8.5]{1}
  \end{scope}
  \begin{scope}[rotate=0]
    \grEmptyCycle[prefix=d,RA=10]{1}
  \end{scope}
  \begin{scope}[rotate=0]
    \grEmptyCycle[prefix=e,RA=11.5]{1}
  \end{scope}
  \begin{scope}[rotate=0]
    \grEmptyCycle[prefix=f,RA=13]{1}
  \end{scope}
  \begin{scope}[rotate=0]
    \grEmptyCycle[prefix=g,RA=14.5]{1}
  \end{scope}
  \begin{scope}[rotate=0]
    \grEmptyCycle[prefix=h,RA=16]{1}
  \end{scope}
  \begin{scope}[rotate=12]
    \grEmptyCycle[prefix=c1,RA=8.5]{1}
  \end{scope}
  \begin{scope}[rotate=10.15]
    \grEmptyCycle[prefix=d1,RA=10]{1}
  \end{scope}
  \begin{scope}[rotate=8.85]
    \grEmptyCycle[prefix=e1,RA=11.5]{1}
  \end{scope}
  \begin{scope}[rotate=7.8]
    \grEmptyCycle[prefix=f1,RA=13]{1}
  \end{scope}
  \begin{scope}[rotate=7]
    \grEmptyCycle[prefix=g1,RA=14.5]{1}
  \end{scope}
  \begin{scope}[rotate=6.35]
    \grEmptyCycle[prefix=h1,RA=16]{1}
  \end{scope}

\node [scale=1.7] at (16,-0.8) {$v_7$};
\node [scale=1.7] at (14.5,-0.8) {$v_6$};
\node [scale=1.7] at (13,-0.8) {$v_5$};
\node [scale=1.7] at (11.5,-0.8) {$v_4$};
\node [scale=1.7] at (10,-0.8) {$v_3$};
\node [scale=1.7] at (8.5,-0.8) {$v_2$};
\node [scale=1.7] at (7,-0.8) {$v_1$};

\node [scale=1.7] at (15.8,2.5) {$u_7$};
\node [scale=1.7] at (14.3,2.5) {$u_6$};
\node [scale=1.7] at (12.8,2.5) {$u_5$};
\node [scale=1.7] at (11.3,2.5) {$u_4$};
\node [scale=1.7] at (9.8,2.5) {$u_3$};
\node [scale=1.7] at (8.3,2.5) {$u_2$};

\draw[dotted, thick](3.8,0)--(6.5,0);
\draw[dotted, thick](2.5,-2.5)--(6.5,0);
\draw[dotted, thick](2.5,2.5)--(6.5,0);

\draw(b0)--(c0)--(d0)--(e0)--(f0)--(g0)--(h0);
\draw(b0)--(c10)--(c0)--(d10)--(d0)--(e10)--(e0)--(f10)--(f0)--(g10)--(g0)--(h10)--(h0);
\end{tikzpicture}
\caption{The graph $H_8$. Edges joining the sets $A$ and $B$, as well as joining $B$ with the vertex $v_1$ are indicated with dotted lines.}\label{fig:H-8}
\end{figure}
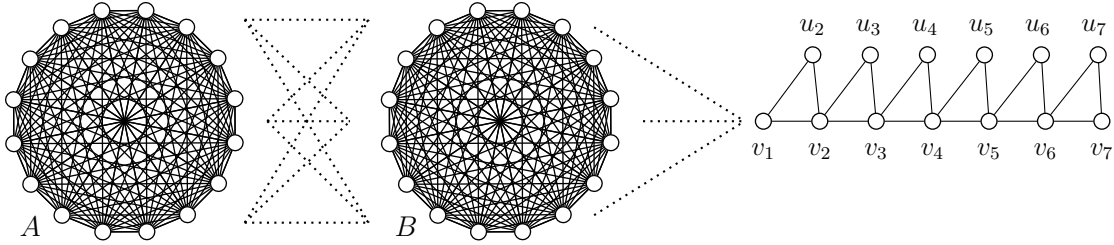

\medskip\noindent
{\bf Strict inequality in~\eqref{eq:chain} in every case}. \\
Given a positive integer $t$, construct the graph $H_t$ as follows. Begin with a complete graph $K_{4t}$ with vertex set $V(K_{4t})=A\cup B$ where $|A|=|B|=2t$. Next, add a path $P_{t-1}=v_1\ldots v_{t-1}$, and join with an edge every vertex of $B$ with the leaf $v_1$ of $P_{t-1}$. Then, add a pendant vertex $u_i$ to every vertex $v_i\in \{v_2,\dots,v_{t-1}\}$, and finally, for every $i\in\{2,\dots,t-1\}$, add the edge $u_iv_{i-1}$. As an example, the graph $H_8$ is represented in Figure \ref{fig:H-8}.

Notice that the graph $H_t$ has diameter $t$. The general $d$-position number of $H_t$ for all possible $d$ is given in the following result.

\begin{proposition}
\label{prop:Ht}
If $2\le d\le t$, then
$$\gp_d(H_t)=\left\{\begin{array}{ll}
                                                    4t;   & d=t, \\
                                                    4t+2; & d=t-1, \\
                                                    5t-d+1; & \mbox{otherwise}.
                                                  \end{array}
\right.$$
\end{proposition}

\proof
We first note that the set $A\cup B$ is a general position set of $H_t$, or equivalently a general $t$-position set. Thus, $\gp(H_t)=\gp_t(H_t)\ge 4t$. Suppose $\gp(H_t)=\gp_t(H_t)>4t$ and let $S$ be a general $t$-position set. Hence, there exists at least one vertex not in $A\cup B$ which is in $S$. Since every shortest path joining a vertex of $A$ with a vertex not in $A\cup B$ passes through a vertex in $B$, it follows that $S\cap A=\emptyset$ or  $S\cap B=\emptyset$. This implies that
 $|S|\le 2t+2t-3=4t-3$, and this is not possible. Therefore  $\gp(H_t)=\gp_t(H_t)=4t=5t-d$, when $d=t$.

We next consider the case $d=t-1$. The set $A\cup B\cup \{v_{t-1},u_{t-1}\}$ is a general $(t-1)$-position set of $H_t$, and so, $\gp_{t-1}(H_t)\ge 4t+2$. If we suppose that $\gp_{t-1}(H_t)> 4t+2$, then a similar argument to that above for $d=t$ leads to a contradiction. Therefore, $\gp_{t-1}(H_t)=4t+2$.

We finally consider $d=t-k$ with $2\le k\le t-2$. Notice that the set $A\cup B\cup \{v_{t-1},u_{t-1},u_{t-2},\dots,u_{t-k}\}$ is a general $d$-position set of $H_t$ of cardinality $4t+k+1=4t+(t-d+1)=5t-d+1$, and so, $\gp_d(H_t)\ge 5t-d+1$. Again, an argument similar to the two cases above leads to $\gp_d(H_t)=5t-d+1$.
\qed

Proposition~\ref{prop:Ht} yields strict inequalities in the chain \eqref{eq:chain}, that is, for any graph $H_t$ with $t\ge 3$, we have
\begin{equation}
\gp_{t}(H_t) < \gp_{t-1}(H_t) < \cdots < \gp_{2}(H_t)\,.
\end{equation}

We shall finish this section by considering the computational complexity of the decision problem related to finding the general $d$-position number of graphs, in which we also show the usefulness of the above graphs $H_t$.

\begin{center}
	\fbox{\parbox{0.80\linewidth}{\noindent
			{\sc General $d$-Position Problem}\\[.8ex]
			\begin{tabular*}{0.95\textwidth}{rl}
				{\em Input:} &  A graph $G$, an integer $d\ge 2$, and a positive integer $r$.  \\
				{\em Question:} & Is $\gp_d(G)$  larger than $r$? \\
			\end{tabular*}
	}}
\end{center}

We first remark that the {\sc General $d$-Position Problem} is known to be NP-complete for every $d\ge \diam(G)$ (see \cite{manuel-2018a}). Hence, we may center our attention on the cases $d\in \{2,\dots,\diam(G)-1\}$, although our reduction also works for the case $d=\diam(G)$.

\begin{theorem}
If $d\ge 2$, then the {\sc General $d$-Position Problem} is NP-complete.
\end{theorem}

\proof
First, we can readily observe that the problem belongs to the class NP, since checking that a given set is indeed a general $d$-position set can be done in polynomial time. From now on, we make a reduction from the {\sc Maximum Clique Problem} to the {\sc General $d$-Position Problem}.

In order to present the reduction, for a given graph $G$ of order $t$, we shall construct a graph $G'$ by using the above graphs $H_t$. We construct $G'$ from the disjoint union of $G$ and $H_t$, by adding all possible edges between $A\cup B\cup \{v_1\}$ and $V(G)$. It is then easily observed that $\omega(G')=|A|+|B|+\omega(G)$.  Moreover, using similar arguments as in the proof of Proposition~\ref{prop:Ht}, we deduce that $\gp_d(G')=\gp_d(H_t)+\omega(G)$. From this fact, since the value $\gp_d(H_t)$ is known from Proposition \ref{prop:Ht}, the reduction is completed, and the theorem is proved.
\qed

\section{Paths and cycles}
\label{sec:paths-and-cycles}

In this section we determine the general $d$-position number of paths and cycles. The first result in turn implies a general upper bound on the general $d$-position number in term of the diameter of a given graph.

\begin{proposition}
\label{prop:paths}
If $n\ge 3$ and $2\le d\le n-1$, then
$$\gp_d(P_n)=\left\{\begin{array}{ll}
2\left\lceil\frac{n}{d+1}\right\rceil-1; & n\equiv 1 \pmod{d+1}, \\[0.2cm]
2\left\lceil\frac{n}{d+1}\right\rceil; & \mbox{otherwise}.
\end{array}
\right.$$
\end{proposition}

\proof
Let $d\in\{2,\dots,n-1\}$ and $P_n=v_1v_2\,\ldots\, v_n$. If $n\equiv 1 \pmod{d+1}$, then let
$$S=\{v_{(d+1)i+1},v_{(d+1)i+2}:\ 0\le i\le \left\lfloor n/(d+1)\right\rfloor-1\}\cup \{v_n\}\,,$$
and if  $n\not \equiv 1 \pmod{d+1}$, then let
$$S=\{v_{(d+1)i+1},v_{(d+1)i+2}:\ 0\le i\le \left\lceil n/(d+1)\right\rceil - 1\}\,.$$
It can be readily seen that  $S$ is a general $d$-position set of $P_n$, which gives the lower bound
$$\gp_d(P_n)\ge \left\{\begin{array}{ll}
2\left\lceil\frac{n}{d+1}\right\rceil-1; & n\equiv 1 \pmod{d+1}, \\[0.2cm]
2\left\lceil\frac{n}{d+1}\right\rceil; & \mbox{otherwise}.
\end{array}
\right.$$
On the other hand, suppose
$$\gp_d(P_n)>\left\{\begin{array}{ll}
2\left\lceil\frac{n}{d+1}\right\rceil-1; & n\equiv 1 \pmod{d+1}, \\[0.2cm]
2\left\lceil\frac{n}{d+1}\right\rceil; & \mbox{otherwise},
\end{array}
\right.\,$$
and let $S'$ be a general $d$-position set of cardinality $\gp_d(P_n)$.  By the pigeonhole principle, we deduce that there exists a subpath in $P_n$ of length $d$ that contains at least three elements of $S'$, but this is not possible. Therefore, the desired equality follows.
\qed

Specializing to $n=14$ in Proposition~\ref{prop:paths}, we next show a table with the values of $\gp_d(P_n)$ for every possible value of $d$. Notice that, equalities and inequalities occur in distinct positions with respect to the chain \eqref{eq:chain}.

\begin{table}[ht!]
  \centering
\begin{tabular}{|c|c|c|c|c|c|c|c|c|c|c|c|c|}
  \hline
  $d$ & 2 & 3 & 4 & 5 & 6 & 7 & 8 & 9 & 10 & 11 & 12 & 13 \\ \hline
  $\gp_d(P_{14})$ & 10 & 8 & 6 & 6 & 4 & 4 & 4 & 4 & 4 & 4 & 3 & 2 \\
  \hline
\end{tabular}
  \caption{The values of $\gp_d(P_{14})$ for every $2\le d\le 13$.}\label{tab:gp_d-P_14}
\end{table}

The result for paths gives the following general lower bound.

\begin{corollary} \label{cor:diameter}
Let $G$ be a connected graph of diameter $d$.  If $2 \le k \le d$, then
$$\gp_k(G) \ge  \left\{\begin{array}{ll}
2\left\lceil\frac{d+1}{k+1}\right\rceil-1; & d\equiv 0 \pmod{k+1}, \\[0.2cm]
2\left\lceil\frac{d+1}{k+1}\right\rceil; & \mbox{otherwise}.
\end{array}
\right.$$
\end{corollary}

\proof
Shortest paths are isometric subgraphs; in particular, this holds for diametrical paths. Hence $G$ contains an isometric $P_{d+1}$, and therefore $\gp_k(G) \ge \gp_k(P_{d+1})$ by Proposition~\ref{prop:isometric} with $r=1$.  Applying Proposition~\ref{prop:paths} yields the result.
\qed

In a similar manner as done for paths, we can compute the general $d$-position number for cycles. It is easy to show that $\gp_d(C_3)=3$ for any $d$, $\gp_1(C_4)=4$, and $\gp_d(C_4)=2$ for $d \ge 2$.

\begin{proposition}
\label{prop:cycles}
If $n \ge 5$ and $2\le d < \left\lfloor\frac{n}{2}\right\rfloor$, then
\begin{equation*}
\gp_d(C_n)=
   \begin{cases}
          2\left\lfloor\frac{n}{d+1}\right\rfloor+1; & n\equiv d \pmod{d+1}, \\ [0.2cm]
          2\left\lfloor\frac{n}{d+1}\right\rfloor; & \mbox{otherwise.}
   \end{cases}
\end{equation*}
If $d \ge \left\lfloor\frac{n}{2}\right\rfloor$, then $\gp_d(C_n)=3$.
\end{proposition}

\proof
Let $C_n=v_1v_2\ldots v_nv_1$. Note that $\diam(C_n)=\left\lfloor\frac{n}{2}\right\rfloor$, and the argument naturally splits into two cases.

First assume that $2\le d < \left\lfloor\frac{n}{2}\right\rfloor$.  Let $m=\left\lfloor\frac{n}{d+1}\right\rfloor$ and for each $k\in [m]$ we define $X_k$ by
$X_k=\{v_i:\, (k-1)(d+1)+1 \le i \le k(d+1)\}$.  Let $X=V(C_n)\setminus \bigcup_{k=1}^mX_k$.  Note that $|X|=x$ where $n \equiv x \pmod{d+1}$ and $x$ is the unique integer such that
 $0 \le x \le d$.  If $x \neq d$, then let $S=\{v_{(k-1)(d+1)+1}, v_{(k-1)(d+1)+2}:\, 1 \le k \le m\}$.  If $x=d$, then let
 $S=\{v_{(k-1)(d+1)+1}, v_{(k-1)(d+1)+2}:\, 1 \le k \le m\} \cup \{v_{m(d+1)+1}\}$.   It is straightforward to check that in both cases $S$ is a general $d$-position set,
 which shows that the claimed value is a lower bound for $\gp_d(C_n)$.  As in the proof of Proposition~\ref{prop:paths}, an application of the pigeonhole principle establishes the upper bound.

Since $\diam(C_n)=\left\lfloor\frac{n}{2}\right\rfloor$, to prove the second statement it is sufficient to show that $\gp_d(C_n)=3$ for $d=\left\lfloor\frac{n}{2}\right\rfloor$.  For this purpose, let $S=\{v_1,v_3,v_{\left\lceil\frac{n}{2}\right\rceil+2}\}$.  For $n=2r$, we see that $S=\{v_1,v_3,v_{r+2}\}$ and
 $d=r$.  On the other hand, for $n=2r+1$, we have $S=\{v_1,v_3,v_{r+3}\}$ and $d=r$.  In both cases an easy computation shows that none of the three vertices lies on a shortest path in $C_n$ between the other two vertices.  Therefore, $S$ is a general $d$-position set, and it follows that $\gp_d(C_n)\ge 3$.  Suppose $T$ is an arbitrary  general  $d$-position set of $C_n$.  We may assume without loss of generality that $v_1 \in T$.  It follows that $|T \cap \{v_2,\ldots,v_{r+1}\}| \le 1$ and  $|T \cap \{v_{r+2}, \ldots v_n\}| \le 1$, for otherwise $T$ contains three vertices that lie on a path of length at most $d$.  Therefore, $\gp_d(C_n)\le |T| \le 3$. \qed

\section{A characterization of general $d$-position sets}
\label{sec:characterization}

In~\cite[Theorem 3.1]{anand-2019} a structural characterization of general position sets of a given graph was proved. In this section we give such a characterization for general $d$-position sets and as a consequence deduce the characterization from~\cite{anand-2019}.

\begin{theorem}
\label{thm:characterize-d-general-sets}
Let $G$ be a connected graph and let $d\ge 2$ be an integer. Then $S\subseteq V(G)$ is a general $d$-position set if and only if the following conditions hold:
\begin{enumerate}
\item[(i)] $\langle S\rangle$ is a disjoint union of complete graphs $Q_1, \ldots, Q_\ell$.
\item[(ii)] If $Q_i\nparallel Q_j$, $i\ne j$, then $d_G(Q_i, Q_j) \ge d$.
\item[(iii)] If $d_G(Q_i, Q_j) + d_G(Q_j, Q_k) = d_G(Q_i, Q_k)$ for $\{i,j,k\}\in \binom{[\ell]}{3}$, then $d_G(Q_i, Q_k) > d$.
\end{enumerate}
\end{theorem}

\proof
Let $S$ be a general $d$-position set of $G$ and let $H$ be a connected component of $\langle S\rangle$. If $H$ is not complete, then it contains an induced $P_3$. The vertices of this $P_3$ are on a geodesic of length $2$ which is not possible since they belong to $S$ and $d\ge 2$. Hence $H$ must be complete.

Consider next two cliques $Q_i$ and $Q_j$ that are not parallel. Let $d_G(Q_i,Q_j) = p$ and let $u\in Q_i$ and $v\in Q_j$ be vertices with $d_G(u,v) = p$. Since $Q_i\nparallel Q_j$, we may assume without loss of generality that there is a vertex $w\in Q_i$ such that $d_G(w,Q_j) = p+1$. Then $u$ lies on a $w,v$-geodesic of length $p+1$ which implies that $p+1 \ge d+1$ and so, $d_G(Q_i,Q_j) = p \ge d$.

Assume next that $d_G(Q_i, Q_j) + d_G(Q_j, Q_k) = d_G(Q_i, Q_k)$ for some $\{i,j,k\}\in \binom{[\ell]}{3}$. If $Q_i\nparallel Q_j$, then by the already proved condition (ii) we immediately get that $d_G(Q_i, Q_j) \ge d$ and thus $d_G(Q_i, Q_k) > d$. The same holds if $Q_j\nparallel Q_k$. Hence assume next that $Q_i\parallel Q_j$ and $Q_j\parallel Q_k$. Let $u\in Q_i$ and $w\in Q_k$ be vertices with $d_G(u,w) = d_G(Q_i,Q_k)$. Since  $d_G(Q_i, Q_j) + d_G(Q_j, Q_k) = d_G(Q_i, Q_k)$, $Q_i\parallel Q_j$, and $Q_j\parallel Q_k$, it follows that  $d_G(u,w) = d_G(u,v) + d_G(v,w)$  for every vertex $v$ of $Q_j$. We conclude that $d_G(Q_i, Q_k) > d$.

\medskip
To prove the converse, assume that conditions (i), (ii), and (iii) are fulfilled for a given set $S$ and let $\{u,v,w\}\in \binom{S}{3}$. We need to show that $u,v,w$ lie in general $d$-position.

If $u,v,w$ lie in the same connected component of $\langle S\rangle$, then by (i), this component is complete and the assertion is clear. Suppose next that $u,v,w$ lie in
the union of cliques $Q_i$ and $Q_j$. If $Q_i\parallel Q_j$, then $u,v,w$ are clearly in general $d$-position. And if $Q_i\nparallel Q_j$, then $u,v,w$ lie in general $d$-position by (ii).

In the last case to consider the three vertices lie in different cliques, say $u\in Q_i$, $v\in Q_j$, and $w\in Q_k$. If the assertion does not hold, then the three vertices lie on a common geodesic and we may assume without loss of generality that $d_G(u,w) = d_G(u,v) + d_G(v,w)$. If $Q_i\nparallel Q_j$, then by (ii), we get $d_G(Q_i, Q_j) \ge d$ and hence $d_G(u,w) = d_G(u,v) + d_G(v,w) \ge d_G(Q_i,Q_j) + d_G(Q_j,Q_k) \ge d + 1 > d$. Analogously, if $Q_j\nparallel Q_k$, we also get $d_G(u,w) > d$. Suppose then that $Q_i\parallel Q_j$ and $Q_j\parallel Q_k$. If also $Q_i\parallel Q_k$, then $d_G(u,w) = d_G(u,v) + d_G(v,w)$ implies that $d_G(Q_i, Q_j) + d_G(Q_j, Q_k) = d_G(Q_i, Q_k)$ and so $d_G(Q_i, Q_k) > d$ by (iii). Again using the fact that $Q_i\parallel Q_k$, it follows that $d_G(u,w) > d$. We are left with the case that $Q_i\parallel Q_j$, $Q_j\parallel Q_k$, and $Q_i\nparallel Q_k$. If $d_G(u,w) = d_G(Q_i,Q_k)$, then by (iii), we get that $d_G(u,w) > d$. Otherwise we may assume without loss of generality that there exists a vertex $u'\in Q_i$, $u'\ne u$, such that $d_G(Q_i,Q_k) = d_G(u',Q_k) < d_G(u,w)$. Since $Q_i\nparallel Q_k$, (ii) implies that $d_G(u',Q_k)\ge d$. But then $d_G(u,w) > d_G(Q_i,Q_k) \ge d$.
\qed

\begin{corollary} {\rm \cite[Theorem 3.1]{anand-2019}}
Let $G$ be a connected graph. Then $S\subseteq V(G)$ is a general position set if and only if the following conditions hold:
\begin{enumerate}
\item[(i)] $\langle S\rangle$ is a disjoint union of complete graphs $Q_1, \ldots, Q_\ell$.
\item[(ii)] $Q_i\parallel Q_j$ for every $i\ne j$.
\item[(iii)] $d_G(Q_i, Q_j) + d_G(Q_j, Q_k) \ne d_G(Q_i, Q_k)$ for every $\{i,j,k\}\in \binom{[\ell]}{3}$.
\end{enumerate}
\end{corollary}

\proof
Set $d=\diam(G)$, so that general $d$-position sets are precisely general position sets. Condition (ii) of Theorem~\ref{thm:characterize-d-general-sets} implies that in cliques $Q_i$ and $Q_j$, which are not parallel, we can find a pair of vertices at distance larger than $\diam(G)$.  Since this is not possible, every two cliques must be parallel. Similarly, if the assumption of condition (iii) would be fulfilled for some cliques $Q_i$, $Q_j$, and $Q_k$, then we would again have vertices at distance larger than $\diam(G)$. Therefore, $d_G(Q_i, Q_j) + d_G(Q_j, Q_k) \ne d_G(Q_i, Q_k)$ must hold for every $\{i,j,k\}\in \binom{[\ell]}{3}$.
\qed

\section{Connections with other topics}
\label{sec:connections}

In this section we connect general $d$-position sets with the dissociation number and with strong resolving graphs.

\subsection*{Strong resolving graphs}

A vertex $u$ of a connected graph $G$ is \emph{maximally distant} from a vertex $v$ if every $w\in N(u)$ satisfies $d_G(v,w)\le d_G(u,v)$. If $u$ is maximally distant from $v$, and $v$ is maximally distant from $u$, then $u$ and $v$ are \emph{mutually maximally distant} (MMD for short).  Given an integer $d\ge 2$, the {\em strong $d$-resolving graph} $G\sr^d$ of $G$ has vertex set $V(G)$, and two vertices $u,v$ are adjacent in $G\sr^d$ if either $u,v$ are MMD in $G$, or $d_G(u,v)\ge d$. The terminology used in this construction comes from the notion of the strong resolving graph introduced in~\cite{Oellermann2007} as a tool to study the strong metric dimension of graphs.  See also~\cite{Kuziak-2018}.

The following observation will be useful in the proof of Theorem~\ref{thm:lower-bound-SRG}.
\begin{observation} \label{obs:MD-no-geodesic}
If $G$ is connected and a vertex $u$ of $G$ is maximally distant from a vertex $v$ of $G$, then $u \notin I(v,w)$ for every  $w\in V(G)\setminus \{u\}$.
\end{observation}

\proof
For the sake of contradiction suppose there exists such a vertex $w\in V(G)\setminus \{u\}$ such that $u \in I(v,w)$.  Suppose that $v=v_0\ldots v_{i-1}u=v_iv_{i+1}\ldots v_k=w$ is a $v,w$-geodesic.  Since this is a geodesic, it follows that $d(v,u)=i$.  But $u$ is maximally distant from $v$, and thus $d(v,v_{i+1}) \le d(v,u)=i$.  Now, by following a shortest $v,v_{i+1}$-path with the path $v_{i+2} \ldots v_k=w$ we arrive at a $v,w$-path of length less than $k$, which is a contradiction.
\qed

From Observation~\ref{obs:MD-no-geodesic} it follows immediately that if three vertices $x,y,z$ are pairwise  MMD, then $x \notin I(y,z)$, $y \notin I(x,z)$, and $z \notin I(x,y)$.  From this we infer that $x,y,z$ lie in general $d$-position.

\begin{theorem}
\label{thm:lower-bound-SRG}
If $G$ is a connected graph and $d\ge 2$ is an integer, then $\gp_d(G)\ge \omega(G\sr^d)$.
\end{theorem}

\proof
We consider a set $S\subseteq V(G\sr^d)$ that induces a (largest) complete subgraph of $G\sr^d$. Then every two vertices $x,y\in S$ are MMD in $G$, or $d_G(x,y)\ge d$. We now consider three vertices $x,y,z$ of $S$ in the graph $G$. If they are pairwise MMD in $G$, then as above, $x,y,z$ lie  in general $d$-position.  Suppose then that two of them, say $x$ and $y$, are not MMD in $G$.  Since $x,y$ are adjacent in $G\sr^d$, it follows that $d_G(x,y)\ge d$.  Suppose for instance that $x,z$ are MMD in $G$.  By Observation~\ref{obs:MD-no-geodesic}, it follows that $x \notin I(z,y)$ and $z \notin I(x,y)$.   If $y\in I(x,z)$, then $d_G(x,z)=d_G(x,y)+d_G(y,z)\ge d+1$, and hence $x,y,z$ lie in general $d$-position.  On the other hand, if $y\notin I(x,z)$, then by definition, $x,y,z$ lie in general $d$-position.

It remains only to consider the case in which no pair of $x,y,z$ is MMD in $G$. This means that the distance between any two of them is at least $d$, and this clearly means that $x,y,z$ are in general $d$-position.
\qed

Note that if $d=\diam(G)$, then $G\sr^d$ is the standard strong resolving graph $G\sr$ as defined in~\cite{Oellermann2007}. In this case Theorem~\ref{thm:lower-bound-SRG} reduces to $\gp(G)\ge \omega(G\sr)$, a result earlier obtained in~\cite[Theorem 3.1]{klavzar-2019}.

\subsection*{Dissociation number and independence number}

If $G$ is a graph and $S\subseteq V(G)$, then $S$ is a {\em dissociation set} if $\langle S\rangle$ has maximum degree at most $1$. The {\em dissociation number} $\diss(G)$ of $G$ is the cardinality of a largest dissociation set in $G$. This concept was introduced by Yanakkakis~\cite{yannakakis-1981}; see also~\cite{boliac-2004, bresar-2017, kardos-2011}. Further, a {\em $k$-path vertex cover} of $G$ is a subset $S$ of vertices of $G$ such that every path of order $k$ in $G$ contains at least one vertex from $S$. The minimum cardinality of a $k$-path vertex cover in $G$ is denoted by $\psi_k(G)$. The minimum $3$-path vertex cover is a dual problem to the dissociation number because $\diss(G) = |V(G)| - \psi_3(G)$; see~\cite{kardos-2011}. For the algorithmic state of the art on the $3$-path vertex cover problem see~\cite{bai-2019}.

\begin{proposition}
\label{prop:dissociation-gp_2}
If $G$ is a triangle-free graph, $d\ge 2$, and $S\subseteq V(G)$ is a general $d$-position set, then $S$ is a dissociation set. Moreover, if $d=2$, then $S$ is a general $2$-position set if and only if $S$ is a dissociation set.
\end{proposition}

\proof
Let $d$ be a positive integer such that $d \ge 2$.  Suppose that $S$ is a general $d$-position set in a triangle-free graph $G$.  By Theorem~\ref{thm:characterize-d-general-sets} every component of the subgraph $\langle S\rangle$ of $G$ induced by $S$ is a complete graph.  Since $G$ is triangle-free, we conclude that each of these components has order $1$ or $2$.  Therefore, $S$ is a dissociation set.  Now assume that $d=2$ and $S$ is a dissociation set in $G$.  The components $C_1, \ldots, C_k$ of $\langle S\rangle$ each have order $1$ or $2$ and are thus complete graphs.  For every pair of distinct indices $i,j$ in $[k]$, the fact that $C_i$ and $C_j$ are distinct components of the induced subgraph $\langle S\rangle$ implies that $d_G(C_i,C_j) \ge 2$.  Therefore, conditions $(ii)$ and $(iii)$ of Theorem~\ref{thm:characterize-d-general-sets} follow immediately, and hence $S$ is a general $2$-position set.
\qed

Proposition~\ref{prop:dissociation-gp_2} immediately gives the following result for triangle-free graphs.

\begin{corollary}
\label{cor:diss-gp_2-number}
If $G$ is a triangle-free graph and $d\ge 2$, then $\gp_d(G) \le \diss(G)$. Moreover, $\gp_2(G) = \diss(G)$.
\end{corollary}

We next relate the particular case of general $2$-position number with the independence number of graphs.

\begin{proposition}
\label{prop:no-true-twins}
If $G$ is a connected graph without true twins, then $\gp_2(G)\ge \alpha(G)$.
\end{proposition}

\proof
Let $x,y \in V(G)$. Suppose first that $xy \in E(G)$.  Since $x$ and $y$ are not true twins, it follows that $x$ and $y$ are not MMD.  By definition, we infer that $xy \notin E(G\sr^2)$.  On the other hand, if $xy \notin E(G)$, then $d_G(x,y) \ge 2$ and by definition $xy \in E(G\sr^2)$.  Consequently, $G\sr^2$ is the complement $\overline{G}$ of $G$. Then by using Theorem~\ref{thm:lower-bound-SRG}, we have $\gp_2(G)\ge\omega(\overline{G})=\alpha(G)$.
\qed

It is straighforward to see that if $2\le m\le n$, then $\gp_{2}(K_{m,n}) = n = \alpha(K_{m,n})$. Hence the bound of Proposition~\ref{prop:no-true-twins} is sharp. For another such family consider the grid graphs $P_{2r}\cp P_{2s}$. (For the definition of the Cartesian product operation $\cp$ see, for instance,~\cite{imklra-2008}.)
As already mentioned,  $\psi_3(G)=n-\diss(G)$ holds for any graph $G$ of order $n$. Also, from \cite{bresar-2013} it is known that $\psi_3(P_{2r}\cp P_{2s})=2rs$. Moreover, from Corollary~\ref{cor:diss-gp_2-number}, we have that $\gp_2(P_{2r}\cp P_{2s})=\diss(P_{2r}\cp P_{2s})$. Thus,
$$\gp_2(P_{2r}\cp P_{2s})=\diss(P_{2r}\cp P_{2s})=4rs-\psi_3(P_{2r}\cp P_{2s})=2rs=\alpha(P_{2r}\cp P_{2s})\,.$$

\section{Infinite graphs and some open problems}
\label{sec:conclude}

The general position problem has been partially studied also on infinite graphs. In~\cite{manuel-2018b} it was proved that   ${\rm gp}({P_\infty^2}) = 4$, where $P_\infty^2$ is the $2$-dimensional grid graph (alias the Cartesian product of two copies of the two way infinite path). The general position number of the $2$-dimensional strong grid graph was also determined, and it was shown that $10\le {\rm gp}({P_\infty^3}) \le 16$. In~\cite{klavzar-2019+} the latter lower bound was improved to $14$. All these efforts were recently rounded off in~\cite{klavzar-2020+} where it is proved that if $n\in {\mathbb N}$, then $\gp(P_\infty^n) = 2^{2^{n-1}}$. On the other hand, the following result reduces the study of the general $d$-position number of infinite graphs to the case $d = \infty$.

\begin{proposition}
If $G$ is an infinite graph and $d < \infty$, then $\gp_d(G) = \infty$.
\end{proposition}

\proof
Let $d < \infty$ be a fixed positive integer. There is nothing to be proved if $d=1$, hence assume that $d\ge 2$.

Suppose first that $\diam(G) = \infty$. In this case $G$ contains an infinite isometric path $P = v_1v_2\ldots$.  It is clear that $\{v_{di}: \ i\in {\mathbb N} \}$ is a general $d$-position set, and hence  $\gp_d(G) = \infty$.

Suppose second that $\diam(G) < \infty$.  Considering an arbitrary vertex of $G$ and its distance levels we infer that $G$ contains a vertex $x$ with $\deg(x) = \infty$. Let $H = \langle N[x]\rangle$.  Since $H$ is an infinite graph, Erd\H{o}s-Dushnik-Miller theorem~\cite{dushnik-1941} implies that $H$ contains a (countably) infinite independent set $I$ or an infinite clique $Q$ (of the same cardinality as $H$). If $H$ contains $Q$, then $Q$ is also a clique of $G$, and hence $G$ contains an infinite general $d$-position set. On the other hand, if $H$ contains $I$, then $I$ is also an independent set of $G$. Moreover, having in mind that  $H = \langle N[x]\rangle$, we infer that each pair of vertices of $I$ is at distance $2$ in $G$. This fact in turn implies that $I$ is an infinite general $d$-position set of $G$. We conclude that $\gp_d(G) = \infty$.
\qed

\subsection{Open questions}

In this section  we point out several questions that, in our opinion, are worthy of consideration.

\begin{itemize}

\item In~\cite[Lemma 5.1]{papadimitriou-1982} there is a polynomial algorithm for the dissociation number of trees $T$  and hence for $\gp_2(T)$. On the other hand, $\gp_{\diam(T)}(T)$ can also be efficiently computed. Hence, is it possible to compute in polynomial time $\gp_d(T)$ for any $2< d < \diam(T)$? More generally, what can be done for the case of block graphs? We know that the simplicial vertices of a block graph form a gp-set. Can the algorithm of Papadimitriou and Yannakakis be modified for block graphs?

\item Compare $\diss(G)$ with $\gp_2(G)$ for graphs $G$ with $\omega(G)\ge 3$. Our guess is that these invariants are incomparable in  such graphs. Is there some relationship when $G$ is a block graph?

\item What is $\gp_d(G)$ whenever $G$ is a grid-like graph?  Note that by applying Corollary~\ref{cor:diss-gp_2-number} together with Theorem 4.1 in~\cite{bresar-2013}, one can
find the value of $\gp_2(P_n\cp P_m)$ for any $n$ and $m$.  Find $\gp_d(P_n\cp P_m)$ for $d \ge 3$.  Find the general $d$-position number of a partial grid graph for $d \ge 2$.
\end{itemize}

\section*{Acknowledgements}

We acknowledge the financial support from the Slovenian Research Agency (research core funding No.\ P1-0297 and projects J1-9109, J1-1693,  N1-0095, N1-0108). This research was initiated while the third author was visiting the University of Ljubljana, Slovenia, supported by ``Ministerio de Educaci\'on, Cultura y Deporte'', Spain, under the ``Jos\'e Castillejo'' program for young researchers (reference number: CAS18/00030).

\end{document}